\documentclass[a4paper,12pt]{article}

\usepackage{amssymb}
\usepackage{amsmath}
\usepackage{amsfonts}                                                           
\usepackage{mathrsfs} 
 
\newcommand{\NN}{\mathbb{N}}                                                    
\newcommand{\ZZ}{\mathbb{Z}}

\newcommand{\Biv}{\mathrm{Biv}}
\newcommand{\BFra}{\mathrm{Biv}^{\star}}
\newcommand{\BFSep}{\mathrm{Sep}^{\star}}

\DeclareFontFamily{OT1}{pzc}{}
\DeclareFontShape{OT1}{pzc}{m}{it}{<-> s * [1.200] pzcmi7t}{}
\DeclareMathAlphabet{\mathpzc}{OT1}{pzc}{m}{it} 

\newcommand{\Icg}[2]{\mathrm{ICG}({#2},{#1})}   
\newcommand{\Ene}[2]{\mathpzc{E}({#2},{#1})}   
  
\newcommand{\Emin}[1]{\mathpzc{E}_{\mathrm{min}}({#1})}   
\newcommand{\Emax}[1]{\mathpzc{E}_{\mathrm{max}}({#1})}

\def \C{\hbox{\sf \rlap{\kern.25em \vrule width.1em height1.6ex depth-.1ex}C}}
\def \D{{\sf I\kern-1.5ptD\,}}
\def \K{{\sf I\kern-1.5ptK\,}}
\def \N{{\sf I\kern-1.5ptN\,}}
\def \P{{\sf I\kern-1.5ptP\,}}
\def \Q{\hbox{\sf \rlap{\kern.25em \vrule width.1em height1.6ex depth-.1ex}Q}}
\def \R{{\sf I\kern-1.5ptR\,}}
\def \T{{\sf T\kern-6.5ptT\,}}
\def \Z{{\sf Z\kern-5.0ptZ\,}}

\def\begitem#1 {\bigskip\pagebreak[1]%
     \refstepcounter{subsection}{\nopagebreak[4]%
     \thesubsection\hskip 0.5truecm}
     {\sc#1}\hskip 1pt.\nopagebreak[4]\par\nopagebreak[4]%
      \begin{enumerate}\rm\nopagebreak[4]}
\def\BEGITEMK#1 #2{\bigskip\pagebreak[1]%
      \refstepcounter{subsection}{\nopagebreak[4]%
     \thesubsection\hskip 0.5truecm}\nopagebreak[4]
     {\bf#1}\hskip 1pt.\nopagebreak[4]\par\nopagebreak[4]%
     \medskip\nopagebreak[4]\rm#2\nopagebreak[4]%
     \begin{enumerate}\nopagebreak[4]\rm}
\def\enditem{\end{enumerate}}

\newcommand{\qed}{\hfill \ensuremath{\Box}\bigskip}

\newtheorem{theorem}{Theorem}[section]

\newtheorem{proposition}{Proposition}[section]

\newtheorem{conjecture}{Conjecture}[section]
\newtheorem{example}{Example}[section]

\begin{document}


\title{The maximal energy of classes of integral circulant graphs} 
\author{J.W.~Sander und T.~Sander}
\maketitle

\begin{abstract}

The energy of a graph is the sum of the moduli of the eigenvalues of its adjacency matrix. We study the energy of integral circulant graphs, also called gcd graphs,
which can be characterized by their vertex count $n$ and a set $\cal D$ of divisors of $n$ in such a way that they have vertex set $\mathbb{Z}_n$ and edge set $\{\{a,b\}:\, a,b\in\mathbb{Z}_n,\, \gcd(a-b,n)\in {\cal D}\}$. For a fixed prime power $n=p^s$ and a fixed divisor set size $\vert {\cal D}\vert =r$, 
we analyze the maximal energy among all matching integral circulant graphs. 
Let $p^{a_1} < p^{a_2} < \ldots < p^{a_r}$ be the elements of ${\cal D}$. 
It turns out that the differences $d_i=a_{i+1}-a_{i}$ between the exponents of an
energy maximal divisor set must satisfy certain balance conditions: 
(i) either all $d_i$ equal $q:=\frac{s-1}{r-1}$, or at most the two differences $[q]$ and $[q+1]$ may occur; 
(ii) there are rules governing the sequence $d_1,\ldots,d_{r-1}$ of consecutive differences. 
For particular choices of $s$ and $r$ these conditions already guarantee maximal energy and its value can be computed explicitly.
\end{abstract}

\\

\\

%
%

%
%

\section{Introduction}

Integral circulant graphs have attacted much research attention lately, in particular since
more and more people have become aware that they 
play a role in quantum physics \cite{SAX}, \cite{BAS4}. A characteristic property of circulant graphs
is that their vertices can be numbered such that any cyclic rotation of the vertex numbering results
in a graph isomorphic to the original graph. Circulant graphs have been the object of research
for quite some time \cite{DAV} and belong to the important family of Cayley graphs.
The integral circulant graphs, having only integer eigenvalues, 
form a small but rather distinguished subclass since integral graphs are quite rare among graphs in general \cite{AHM}.

Given an integer $n$ and a set ${\cal D}$ of positive divisors of $n$, the integral circulant graph $\Icg{\cal D}{n}$ is  
defined as the graph having vertex set $\ZZ_n=\{0,1,\ldots,n-1\}$ and edge set $\{ \{a,b\}:~a,b\in \ZZ_n,~ \gcd(a-b,n)\in {\cal D}\}$. 
We consider only loopless gcd graphs, i.e. $n\notin {\cal D}$.
For $\vert {\cal D}\vert =1$ we obtain the subclass of so-called unitary Cayley graphs.
Over the years, the general structural properties of integral circulant graphs have been well 
researched \cite{DEJ}, \cite{BER}, \cite{SO}, \cite{KLO}, \cite{AKH}, 
\cite{BAS}, \cite{KLO2}, \cite{DRO}, \cite{BEA}, \cite{BAS2}.
Due to the connection with quantum physics, emphasis has lately been placed
on researching the energy of integral circulant graphs 
\cite{SHP}, \cite{ILI}, \cite{RAM}, \cite{BAS3}, \cite{PET}, \cite{SA1}, \cite{SA2}.

The \textit{energy} $E(G)$ of a graph $G$ on $n$ vertices is defined as
\[
     E(G)=\sum_{i=1}^n \vert\lambda_i\vert,
\]
 where $\lambda_1,\ldots,\lambda_n$ are the
eigenvalues of the adjacency matrix of $G$.
Refer to \cite{BRU} and \cite{GUT} for general results on graph energy.

Let us abbreviate $\Ene{\cal D}{n}=E(\Icg{\cal D}{n})$.
Given a positive integer $n$, we consider
\[ \Emin{n} := \min\,\{ \Ene{\cal D}{n}:\;\; {\cal D}\subseteq \{1\le d<n:\; d\mid n\} \}\]
and
\[ \Emax{n} := \max\,\{ \Ene{\cal D}{n}:\;\; {\cal D}\subseteq \{1\le d<n:\; d\mid n\} \}.\]

Consider a prime power $n=p^s$ and a divisor set ${\cal D} = \{ p^{a_1}, p^{a_2} , \ldots , p^{a_r}\}$ with exponents $0\le a_1 < \ldots < a_r \le s-1$. 
According to
Theorem 2.1 in \cite{SA1} we have
\begin{equation}
\Ene{{\cal D}}{p^s} = 2(p-1)p^{s-1}\left(r -(p-1)h_p(a_1,\ldots,a_r) \right),  \label{ft3}
\end{equation}
where
\begin{equation}
h_p(x) = h_p(x_1,\ldots,x_r) := \sum_{k=1}^{r-1}\sum_{i=k+1}^r \,\frac{1}{p^{x_i-x_k}}   \label{basishp}
\end{equation}
for $x=(x_1, \ldots, x_r) \in \mathbb{R}^r$. Observe that $h_p$ has the symmetry property
\begin{equation}
 h_p(s-1-a_r,\ldots,s-1-a_1) = h_p(a_1,\ldots,a_r) \label{symmhp}
 \end{equation} 
for all integral exponents $0\leq a_1 < a_2 < \ldots < a_{r-1} < a_r \leq s-1$.
A straightforward consequence of (\ref{ft3}) is that $\Emin{p^s}$
is attained precisely for the singleton divisor sets ${\cal D}=\{p^t\}$ with $0\le t \le s-1$ (cf. \cite{SA1}, Theorem 3.1).

In \cite{SA2} divisor sets $\cal D$ producing graphs with maximal energy 
$ \Emax{p^s}$ were studied. Equivalently, exponent tuples $(a_1,\ldots,a_r)$ minimizing $h_p$ had to be found. 
By the result cited above, such minimizers satisfy $r\geq 2$, and they obviously must have the entries $a_1=0$ and $a_r=s-1$.  Accordingly, a corresponding $a=(a_1,\ldots,a_r)$ lies in the set
\[ A(s,r) := \{(a_1,\ldots,a_r)\in \mathbb{Z}^r:\; 0=a_1<a_2<\ldots <a_{r-1}<a_r=s-1\},\]
and such an $a$ is called an \emph{admissible} exponent tuple.

Hence the quest for minimizers of $h_p$ is only interesting in case $r\geq 3$, which we shall assume in the sequel. It was shown by use of methods from convex optimization that, for fixed $s$ and $r$, the function $h_p$ becomes almost minimal if only $0=a_1<a_2<\ldots<a_{r-1}<a_r=s-1$ are chosen in nearly equidistant position (\cite{SA2}, Theorem 4.2). Note here that perfect equidistance can only be achieved if $(r-1) \mid (s-1)$ because the $a_i$ are integers.
It is the purpose of this article to use combinatorial instead of analytic arguments in order to refine the earlier approximative results. 

The nearly equidistant positioning just mentioned indicates that the key to maximizing the energy lies 
in considering the successive exponent differences. Hence, for a given $a\in A(s,r)$, we define its \emph{delta vector} as
\[
\delta(a):= (\delta_1(a), \delta_2(a), \ldots, \delta_{r-1}(a)) \in \mathbb{N}^{r-1}
\]
with $\delta_j(a):=a_{j+1}-a_j$ ($1\le j \leq r-1$). Obviously, we have $\sum_{j=1}^{r-1} \delta_j(a) = s-1$. Thus, introducing
\[
D(s,r):= \{(d_1,\ldots,d_{r-1})\in \mathbb{N}^{r-1}: \; \sum_{j=1}^{r-1} d_j = s-1 \}, 
\] 
the function
\[ 
\delta: \left\{ \begin{array}{ccl}
                   A(s,r) &  \longrightarrow & D(s,r) \\
                   (a_1,a_2,\ldots,a_r) & \mapsto & (a_2-a_1,a_3-a_2,\ldots,a_r-a_{r-1}) 
                   \end{array}  
                   \right.
\] 
is 
1--1 
with its inverse 
\[ \quad\quad \quad\quad\;\;\,
\delta^{-1}: \left\{ \begin{array}{ccl}
                   D(s,r) &  \longrightarrow & A(s,r) \\
                   (d_1,d_2,\ldots,d_{r-1}) & \mapsto & (0,d_1,d_1+d_2,\ldots, d_1+d_2+\ldots+d_{r-2,s-1}). 
                   \end{array}  
                   \right.                  
\]
%
The mentioned divisor set structure becomes apparent by restrictions
on the delta vector $\delta(a)$ corresponding to an energy maximal exponent tuple $a$ as follows:\newline 
Firstly, the set $\{\delta_j(a):\; j=1,\ldots,r-1\}$ of differences is either a singleton or has only two elements that are successive positive integers. 
Secondly, the distribution of the differences must satisfy certain balance conditions, in the sense that the differences of the value occuring less often than the other must be distributed somewhat ``evenly'' between the other difference values.

In some cases, these restriction will already characterize the
delta vectors, and consequently the divisor set(s) imposing maximal energy on the corresponding class of integral circulant graphs. In other words, for some fixed $s$ and $r$, we will be able to determine precisely 
\[ \min h_p := \min \{h_p(a): a \in A(s,r)\}   \]
along with all admissible $a$ satisfying $h_p(a)= \min h_p$. 

Open questions and conjectures in the final section disclose our view on how a ``perfect balancing'' process might look in order to determine the admissible $a$ satisfying $h_p(a)= \min h_p$ in all cases.

%
%

\section{Main results}

In what follows, we shall consider $3\le r < s$ to be fixed integers and set $q:=\frac{s-1}{r-1}$. 
Furthermore, $p$ will always be a fixed prime.

\bigskip

If $s\equiv 1 \bmod (r-1)$ or $s\equiv 0 \bmod (r-1)$, we are able to determine all minimizers of $\min h_p$ precisely.
\begin{theorem}{\label{Thm2.1}}
Let $p\geq 3$ be a fixed prime, and let $3\leq r < s$.
Assume that $a\in A(s,r)$ is a minimizer of $h_p$, i.e. $h_p(a) = \min h_p$.
\begin{itemize}
\item[(i)] 
If $(r-1)\mid (s-1)$, i.e. $q=\frac{s-1}{r-1}$ is an integer, then 
$a = \delta^{-1}(q,\ldots,q)$,
and we have 
\[ h_p(a) = \min h_p =  \frac{1}{p^q-1}\left(r-1 - \frac{1}{p^q-1}\left(1-\frac{1}{p^{q(r-1)}}\right)\right). \]
\item[(ii)]
If $(r-1)\mid s$, then 
$a=\delta^{-1}([q],[q+1],\ldots,[q+1])$ or $a=\delta^{-1}([q+1],\ldots,[q+1],[q])$,
and we have
\[
  h_p(a) = \min h_p =  \frac{1}{p^{[q+1]}-1}\left(r-1 + \left(p-1-\frac{1}{p^{[q+1]}-1}\right)\left(1-\frac{1}{p^{[q+1](r-1)}}\right)\right).  
\]
\end{itemize} 
\end{theorem}
Inserting the explicit values of $\min h_p$ into formula (\ref{ft3}), one can easily compute the maximal energies of the corresponding classes of integral circulant graphs.
\bigskip\medskip

Complementing Theorem \ref{Thm2.1}, we have the following 
\begin{theorem}{\label{Thm2.2}}
Let $p\geq 3$ be a fixed prime, and let $3\leq r < s$ be such that $(r-1)\nmid (s-1)$ and $(r-1)\nmid s$. Define the integer $g$ as the least positive residue satisfying $g \equiv s-1 \bmod (r-1)$. Assume that $a\in A(s,r)$ is a minimizer of $h_p$, i.e. $h_p(a) = \min h_p$.
\medskip\newline
For $2g\geq r-1$ and $q_2:=\frac{g}{r-g-2}$ we have:
\begin{itemize}
\item[(i)] If $(r-g-2)\mid g$, then 
\[ a= \delta^{-1}\Big([q],\underbrace{[q+1],\ldots,[q+1]}_{q_2-fold},[q],\underbrace{[q+1],\ldots,[q+1]}_{ q_2-fold},[q],\mbox{\rm etc.},\underbrace{[q+1],\ldots,[q+1]}_{q_2-fold},[q]\Big).
\]
\item[(ii)] If $(r-g-2)\nmid g$, then $d=(d_1,\ldots,d_{r-1}):=\delta(a)$ has the following properties:
\begin{itemize}
\item[$\bullet$] There are exactly $r-g-1$ entries $[q]$, two of which are $d_1=d_{r-1}=[q]$. Moreover, neighboring entries $d_j=d_{j+1}=[q]$ do not occur.
\item[$\bullet$]  The remaining $g$ entries of $\delta(a)$ all equal $[q+1]$ and appear in blocks of length either $[q_2]$ or $[q_2+1]$. More precisely, $\delta(a)$ has exactly $\,e\,$ $[q+1]$-blocks of length $[q_2+1]$ and 
$\,(r-g-2-e)\,$ $[q+1]$-blocks of length $[q_2]$, where $e \equiv g \bmod (r-g-2)$ is the least positive residue.  
\end{itemize}
\end{itemize}
For $2g\leq r-2$ and $q_1:=\frac{r-g-1}{g+1}$ we have:
\begin{itemize}
\item[(iii)] If $(g+1)\mid (r-g-1)$, then 
\[ a= \delta^{-1}\Big(\underbrace{[q],\ldots,[q]}_{q_1-fold},[q+1],\underbrace{[q],\ldots,[q]}_{ q_1-fold},[q+1],\mbox{etc.},\underbrace{[q],\ldots,[q]}_{q_1-fold}\Big).
\]
\item[(iv)]  If $(g+1)\nmid (r-g-1)$, then $d=(d_1,\ldots,d_{r-1}):=\delta(a)$ has the following properties:
\begin{itemize}
\item[$\bullet$]  There are exactly $g$ entries $[q+1]$, but $d_1\neq [q+1]$ and  $d_{r-1}\neq  [q+1]$. Moreover, neighboring entries $d_j=d_{j+1}=[q+1]$ do not occur.
\item[$\bullet$]  The remaining $r-g-1$ entries of $\delta(a)$ all equal $[q]$ and appear in blocks of length either $[q_1]$ or $[q_1+1]$. More precisely, $\delta(a)$ has exactly $\,f\,$ $[q]$-blocks of length $[q_1+1]$ and 
$\,(g+1-f)\,$ $[q]$-blocks of length $[q_1]$, where $f \equiv r-g-1 \bmod (g+1)$  is the least positive residue. 

\end{itemize}
 
\end{itemize}
\end{theorem}
As in Theorem \ref{Thm2.1}, the computation of $\min h_p$ in (i) and (iii) is just a matter of evaluating certain multi-geometric sums, and again by use of (\ref{ft3})
this would give explicit formulae for the maximal energies of the corresponding classes of integral circulant graphs.

%
%
\section{ Bivalence -- Proof of Theorem \ref{Thm2.1}}

For $d=(d_1,\ldots,d_{r-1}) \in D(s,r)$, let
\begin{align*}
\max d &:= \max\{d_j: \; 1\le j\le r-1\},  \\ 
\min d &:= \min\{d_j: \; 1\le j\le r-1\}.
\end{align*}
By the definition of $D(s,r)$, we clearly have 
\begin{equation}
1 \le \min d \le q \le \max d \le s-r+1.  \label{diffdelta}
\end{equation}

For 
$d\in D(s,r)$ we call $\rho(d):= \max d - \min d$
the \emph{range} of $d$. 
Any vector containing only entries $m$ or $m+1$ for some positive integer $m$ shall be called {\em bivalent}.
Hence, $d$ 
is bivalent if $\rho(d) \leq 1$. 
It is an immediate consequence of (\ref{diffdelta}) that 
$\min d = [q]$ and $\max d = [q+1]$ for a bivalent $d\in D(s,r)$ in case $q$ is not integral, and $\min d =\max d = q$ for a bivalent $d\in D(s,r)$ in case $q$ is an integer.

%
For the set 
\[ 
\Biv(s,r):= \{d\in D(s,r):\; \rho(d)\le 1\} \subseteq D(s,r),
\] 
containing all bivalent elements of $D(s,r)$, we thus have
\begin{equation} 
\Biv(s,r) = \left\{\begin{array}{ll}
            \{d\in D(s,r):\; \forall j\;\; d_j=[q] \mbox{ or } d_j=[q+1]\}& \mbox{ if $q\notin\mathbb{N}$, } \\
            \{(q,q,\ldots,q)\}& \mbox{ if $q\in\mathbb{N}$. } 
            \end{array} \right.
\label{Abar}
\end{equation}

\begin{proposition}{\label{PropSmallRange}}
Let $p$ be a fixed prime. 
If $a\in A(s,r)$ satisfies $h_p(a) = \min h_p$, 
then $\delta(a) \in \Biv(s,r)$.

\end{proposition}

{\sc Proof. }
We make the assumption that $d=(d_1,\ldots,d_{r-1}):=\delta(a) \notin \Biv(s,r)$ 
and shall derive a contradiction.

Let $u$ be some index such that $d_u = \min d$, and let $v$ be some index such that $d_v= \max d$. By assumption, $u\neq v$. By the symmetry property (\ref{symmhp}) of $h_p$, we may assume w.l.o.g. that $u<v$, and also that 
\begin{equation}
\min d < d_j < \max d   \quad\quad (u<j<v).
\label{intermediate}
\end{equation}

For $a=(a_1,\ldots,a_r)$, say, we define $b=(b_1,\ldots,b_r)\in A(s,r)$ by setting
\begin{equation} \label{bjot}
b_j := \left\{\begin{array}{ll}
                 a_j   & \mbox{for $j\le u$ or $j\ge v+1$,} \\
                 a_j+1 & \mbox{for $u+1\leq j \leq v$,}
                 \end{array} \right. 
\end{equation}
i.e. we simultaneously extend one of the smallest subintervals of the partition $(a_1,\ldots,a_r)$ by $1$ and shorten one of its longest subintervals by $1$, while all other subintervals remain unchanged in length.

Then, by (\ref{basishp}),
\[h_p(a)-h_p(b)= 
\sum_{k=1}^{r-1}\sum_{i=k+1}^r \,\left(\frac{1}{p^{a_i-a_k}} - \frac{1}{p^{b_i-b_k}}\right).\]
According to the definition of $b$ in (\ref{bjot}), 
the two quotients enclosed in parentheses differ from each other only if $1\leq k \leq u$ and $u+1\leq i \leq v$, or if $u+1\leq k \leq v$ and $v+1\leq i \leq r$. Therefore, 
and since $\sum_{k=1}^u p^{a_k-a_u} \geq 1$ and $\sum_{i=v+1}^r p^{a_{v+1}-a_i} < \sum_{j=0}^{\infty}p^{-j} = \frac{p}{p-1}$, 
\begin{align}
h_p(a)&-h_p(b) = \nonumber \\
&= \sum_{k=1}^{u}\sum_{i=u+1}^v \,\left(\frac{1}{p^{a_i-a_k}} - \frac{1}{p^{(a_i+1)-a_k}}\right)
+ \sum_{k=u+1}^{v}\sum_{i=v+1}^r \,\left(\frac{1}{p^{a_i-a_k}} - \frac{1}{p^{a_i-(a_k+1)}}\right) \nonumber\\
&= (p-1)\left(\frac{1}{p} \sum_{k=1}^{u}\, p^{a_k} \sum_{i=u+1}^v \,\frac{1}{p^{a_i}} -
     \sum_{k=u+1}^{v}\,  p^{a_k} \sum_{i=v+1}^r \,\frac{1}{p^{a_i}} \right) \label{diffhp} \\
&=  (p-1)\left(p^{a_u- a_{u+1}-1} \sum_{k=1}^{u}\, \frac{1}{p^{a_u -a_k}} \sum_{i=u+1}^v \,\frac{1}{p^{a_i-a_{u+1}}} - p^{a_v- a_{v+1}} \sum_{k=u+1}^{v}\, \frac{1}{p^{a_v-a_k}} \sum_{i=v+1}^r \,\frac{1}{p^{a_i-a_{v+1}}} \right)  \nonumber\\   
&=  (p-1)\left(p^{-\min d -1} \sum_{k=1}^{u}\, \frac{1}{p^{a_u -a_k}} \sum_{i=u+1}^v \,\frac{1}{p^{a_i-a_{u+1}}}
   - p^{-\max d} \sum_{k=u+1}^{v}\, \frac{1}{p^{a_v-a_k}} \sum_{i=v+1}^r \,\frac{1}{p^{a_i-a_{v+1}}} \right) \nonumber \\  
&> (p-1)\left(p^{-\min d-1} \sum_{i=u+1}^v \,\frac{1}{p^{a_i-a_{u+1}}} 
  - \frac{p^{- \max d +1}}{p-1} \sum_{k=u+1}^{v}\, \frac{1}{p^{a_v-a_k}} \right). \nonumber
\end{align}
Since $\sum_{i=u+1}^v \,p^{a_{u+1}-a_i}\geq 1$ and  $\sum_{k=u+1}^{v}\, p^{a_k-a_v}<\sum_{j=0}^{\infty}p^{-j} = \frac{p}{p-1}$, it follows from (\ref{diffhp}) that
\[
h_p(a)-h_p(b) > (p-1)\left(p^{- \min d -1} - \frac{p^{-\max d +2}}{(p-1)^2}\right).
\]
In case $\rho(d)\geq 3$, i.e. $\min d \leq \max d -3$, we conclude that
\[
h_p(a)-h_p(b) > (p-1) p^{2-\max d} \left( 1 - \frac{1}{(p-1)^2}\right)\geq 0
\]
for all primes $p$, which proves the proposition.

We are left with the case $\rho(d)=2$, i.e. $\min d = \max d -2$. By (\ref{intermediate}), we have 
\[
d_j = \min d +1 = \max d -1 \quad\quad (u<j<v).
\]
Consequently
\[ 
a_i-a_{u+1} = \sum_{j=u+1}^{i-1} d_j = (i-u-1)(\max d -1) \quad\quad (i= u+1,\ldots,v)
\]
and
\[ 
a_v-a_k = \sum_{j=k}^{v-1} d_j  = (v-k)(\max d -1) \quad\quad (k= u+1,\ldots,v).
\]
Hence by (\ref{diffhp})
\begin{align*}
h_p(a)-h_p(b) &> (p-1)p^{1-\max d} \left(\sum_{i=u+1}^v \,\frac{1}{p^{(i-u-1)(\max d-1)}} 
  - \frac{1}{p-1} \sum_{k=u+1}^{v}\, \frac{1}{p^{(v-k)(\max d-1)}} \right)\\
&= (p-1)p^{1-\max d} \sum_{i=u+1}^v \,\frac{1}{p^{(i-u-1)(\max d-1)}} \left(1- \frac{1}{p-1}\right) \geq 0
\end{align*}  
for all primes $p$, which completes our proof.

\qed

{\sc Proof of Theorem \ref{Thm2.1}{\rm(i)}. }\newline 
Let $a\in A(s,r)$ have the property $h_p(a) = \min h_p$. Then we know by Proposition \ref{PropSmallRange}(i) that $\delta(a)\in \Biv(s,r)$. It follows from (\ref{Abar}) that  $\delta(a)=(q,q,\ldots,q)$.

The proof of the formula for $\min h_p$ is an easy exercise with geometric sums.

\qed

\medskip

Up to this point we know that $\min h_p$ can only be attained by admissible tuples $a$ having bivalent delta vectors, 
that is $\delta(a)\in  \Biv(s,r)$. In the sequel, we shall derive further restrictions for minimizers of $h_p$. For $(r-1)\nmid (s-1)$, the number $q$ is not an integer. If $\delta(a)\in \Biv(s,r)$, thus $\delta(a) \in \{[q],[q+1]\}^{r-1}$ by (\ref{Abar}). 

\begin{proposition}{\label{PropUnframed}}
Let $a\in A(s,r)$ satisfy $h_p(a)= \min h_p$, hence $d=(d_1,\ldots,d_{r-1}):=\delta(a)\in \Biv(s,r)$ by Proposition \ref{PropSmallRange}. 
If
$d_1=[q+1]$ or $d_{r-1}=[q+1]$, then
$d=([q],[q+1],[q+1],\ldots,[q+1])$ or
$d=([q+1],[q+1],\ldots,[q+1],[q])$. 
\end{proposition}

{\sc Proof. } 
By the symmetry of $h_p$ (see (\ref{symmhp})), we may assume w.l.o.g. that 
$d_1=[q+1]$.
Clearly, $d_j=[q]$ for at least one $j$. Hence let
\[
d_1 = d_2 = \ldots = d_{\ell} = [q+1], \quad d_{\ell +1} = [q]
\] 
for a suitable $1\leq \ell\leq r-2$. For $a=(a_1,\ldots,a_r)$, say, we define $b=(b_1,\ldots,b_r) \in A(s,r)$ by setting
\[ b_j := \left\{ \begin{array}{cl}
                  a_j & \mbox{ for $j=1$ or $\ell +2 \leq j \leq r$, }\\
                  a_j-1& \mbox{ for $2\leq j \leq \ell +1$. }
                  \end{array} \right.
\]  
Clearly, $\delta(a)\in \Biv(s,r)$ implies $\delta(b)\in \Biv(s,r)$. 

By (\ref{basishp}), we have
\[
h_p(a)-h_p(b) = 
\sum_{k=1}^{r-1}\sum_{i=k+1}^r \,\left(\frac{1}{p^{a_i-a_k}} - \frac{1}{p^{b_i-b_k}}\right).
\]
According to our definition of $b$, 
the two quotients enclosed in parentheses differ from each other only if $k=1$ and $2\leq i \leq \ell+1$, or if $2\leq k \leq \ell+1$ and $\ell+2 \leq i \leq r$. Therefore, 
\begin{align*}
h_p(a)-h_p(b) 
&= \sum_{i=2}^{\ell+1} \,\left(\frac{1}{p^{a_i-a_1}} - \frac{1}{p^{(a_i-1)-a_1}}\right)
 + \sum_{k=2}^{\ell+1}\sum_{i=\ell+2}^r \,\left(\frac{1}{p^{a_i-a_k}} - \frac{1}{p^{a_i-(a_k-1)}}\right) \\
&= (p-1)\left(\frac{1}{p} \sum_{k=2}^{\ell+1}\, p^{a_k} \sum_{i=\ell+2}^r \,\frac{1}{p^{a_i}} -
     p^{a_1} \sum_{i=2}^{\ell+1} \,\frac{1}{p^{a_i}} \right).
\end{align*}
Observe that $a_1=0$ and $a_k= [q+1](k-1)$ for $2\leq k \leq \ell+1$. Hence
\begin{align}\label{hpdiff} 
\begin{split}
h_p(a)-h_p(b) 
&= (p-1)\left(\frac{1}{p} \sum_{k=0}^{\ell-1}\, p^{[q+1](k+1)} \sum_{i=\ell+2}^r \, 
\,\frac{1}{p^{a_i}} 
- \sum_{i=0}^{\ell-1} \,\frac{1}{p^{[q+1](i+1)}} \right)\\
&= \frac{p-1}{p^{[q+1]}-1} 
\left(\frac{p^{[q+1]}}{p} \left(p^{[q+1]\ell}-1\right) \sum_{i=\ell+2}^r 
\,\frac{1}{p^{a_i}} 
- \left( 1- \frac{1}{p^{[q+1]\ell}} \right)\right). 
\end{split}
\end{align}
If $\ell \leq r-3$, we obtain
\[ 
\sum_{i=\ell+2}^r \, \frac{1}{p^{a_i}} > \frac{1}{p^{a_{\ell+2}}} = 
\frac{1}{p^{[q+1]\ell + [q]}}.
\]
Using this lower bound in (\ref{hpdiff}) shows that the righthand side of (\ref{hpdiff}) is positive. Thus $h_p(a) >h_p(b)$, which would contradict the minimality of $h_p(a)$. 
It remains to consider the case $\ell = r-2$, but then $\delta(a)= 
([q+1],[q+1],\ldots,[q+1],[q])$.

\qed

{\sc Proof of Theorem \ref{Thm2.1}{\rm(ii)}. }\newline 
Let $a\in A(s,r)$ satisfy $h_p(a) = \min h_p$. It follows from Proposition \ref{PropSmallRange} that
$(d_1,\ldots,d_{r-1}):=\delta(a)\in \Biv(s,r)$.
The condition $(r-1)\mid s$ means that $s-1\equiv -1 \bmod (r-1)$, hence
$d_j = [q]$ for exactly one $1\leq j \leq r-1$ and $d_j = [q+1]$ otherwise.
By Proposition \ref{PropUnframed}, the condition $h_p(a) = \min h_p$ implies that $d$ 
equals one of the two 
$(r-1)$-tuples given there. Hence $\delta(a)$ has the desired  form.

The proof of the formula for $\min h_p$ is an easy exercise with geometric sums.

\qed

%
%
\section{Separability}

In case $\mbox{$s\equiv 1 \bmod (r-1)$}$ or $s \equiv 0 \bmod (r-1)$, we know all minimizers $a\in A(s,r)$ of $h_p$ by Theorem \ref{Thm2.1}.
If $s$ belongs to another residue class $\bmod\, (r-1)$, we have a further restriction for minimizers of $h_p$. To this end, we shall call any vector {\em framed} if its first and last entry are the same. We indicate that these entries have value $x$, say, by calling the vector \emph{$x$-framed}.

Let
\[ 
\BFra(s,r):= \{ d=(d_1,\ldots,d_{r-1}) \in \Biv(s,r): \; d_1=d_{r-1} = [q]\}.
\]
denote the set of all bivalent, $[q]$-framed delta vectors.

\begin{proposition}{\label{PropFraming}}
Let $p\geq 3$ be a fixed prime, and let
$(r-1)\nmid (s-1)$ and $(r-1)\nmid s$. If $h_p(a) = \min h_p$ for some $a\in A(s,r)$, then 
$\delta(a) \in \BFra(s,r)$.
\end{proposition}

{\sc Proof. }
Proposition \ref{PropSmallRange} tells us that $(d_1,\ldots,d_{r-1}):=\delta(a)\in \Biv(s,r)$. 
The condition
$(r-1)\nmid (s-1)$ implies that $d_j = [q+1]$ for at least one $j$, and $(r-1)\nmid s$ guarantees indices $j_1\neq j_2$ such that $d_{j_1}=d_{j_2}=[q]$. All this shows that $\delta(a)$ cannot be one of the two $(r-1)$-tuples in Proposition \ref{PropUnframed}, which under our  minimality assumption for $h_p(a)$ yields $d_1=d_{r-1} = [q]$.

\qed

For $d\in \BFra(s,r)$, we have $d= \big([q], d_2,\ldots, d_{r-2}, [q] \big)$, where $\mbox{$d_j \in \{[q], [q+1]\}$}$ for all $j$ by (\ref{Abar}). Now we study the sequences of successive $d_j$ of equal value. For suitable positive integers $t_i=t_i(d)$  $(1\leq i \leq 2w+1)$, say, we have
\[ d = 
\Big(\underbrace{[q],\ldots,[q]}_{t_1-fold},\underbrace{[q+1],\ldots,[q+1]}_{t_2-fold},\underbrace{[q],\ldots,[q]}_{t_3-fold},{\rm etc.}, \underbrace{[q+1],\ldots,[q+1]}_{t_{2w}-fold},\underbrace{[q],\ldots,[q]}_{t_{2w+1}-fold}\Big).
\]
To put it another way, $d$ 
is composed of a $[q]$-block of length $t_1$ followed by a $[q+1]$-block of lengths $t_2$ and then alternately by $[q]$-blocks and $[q+1]$-blocks of respective lengths.
Setting $T_{\ell} = T_{\ell}(d) := \sum_{i=1}^{\ell} t_i$ for $1\leq \ell \leq 2w+1$, we have $T_{2w+1}=r-1$ and
\begin{equation} \label{deftj}
\left.
\begin{array}{ccccccccc}
d_1        &=& d_2       &=& \ldots &=& d_{T_1}      &=& [q]  \\
d_{T_1+1}  &=& d_{T_1+2}  &=& \ldots &=& d_{T_2}   &=& [q+1]  \\
d_{T_2+1} &=& d_{T_2+2} &=& \ldots &=& d_{T_3} &=& [q] \\
   \vdots & & \vdots &  & && \vdots  &    &  \vdots\\
d_{T_{2w-1}+1} &=& d_{T_{2w-1}+2} &=& \ldots &=& d_{T_{2w}} &=& [q+1]  \\
d_{T_{2w}+1} &=& d_{T_{2w}+2} &=& \ldots &=& d_{T_{2w+1}} &=& [q]
\end{array} 
\right\} 
\end{equation}

Denote by $g$ the least non-negative integer satisfying $g \equiv s-1 \bmod (r-1)$.  
It is easily seen that
\[g = \# \{1\leq j \leq r-1: \; d_j = [q+1] \}.
\]
In particular, $g$ does not depend on $d$. 
The definition of the $t_i(d)$ clearly implies
\begin{equation}
 \sum_{\ell=0}^w t_{2\ell+1} = r-g-1  \quad \mbox{ and } \quad \sum_{\ell=1}^w t_{2\ell} = g.
\label{sumstl}
\end{equation}
In case $(r-1)\nmid (s-1)$, i.e. $q$ is not integral, we define for $d\in \BFra(s,r)$  the maximal and minimal lengths of $[q]$-blocks and $[q+1]$-blocks, respectively, occurring in  $d$, namely
\[
\begin{array}{rclrcl}
\eta_{\max}(d) &:=& \max\{t_{2\ell+1}:\; 0\le \ell \leq w\},
 &\eta_{\min}(d) &:=& \min\{t_{2\ell+1}:\; 0\le \ell \leq w\}, \\
\theta_{\max}(d) &:=& \max\{t_{2\ell}:\; 1\le \ell \leq w\},
&\theta_{\min}(d) &:=& \min\{t_{2\ell}:\; 1\le \ell \leq w\}.
\end{array}
\]
Then $q_1:= \frac{r-g-1}{w+1}$ and $q_2:=\frac{g}{w}$ are the average lengths of $t_{2\ell+1}$ and $t_{2\ell}$ respectively, i.e. the average lengths of the $[q]$-blocks and the $[q+1]$-blocks, and we obviously have
\begin{equation} \label{eta}
1\le \eta_{\min}(d) \leq q_1 \leq \eta_{\max}(d)
\end{equation}
and
\begin{equation}  \label{theta}
1 \leq \theta_{\min}(d) \leq q_2 \leq \theta_{\max}(d). 
\end{equation}

A bivalent 
vector, containing both entries $m$ and $m+1$, say, shall be called {\em separable} if no consecutive entries $m$ or no consecutive entries $m+1$ occur.
Our next result shows that for an $a\in A(s,r)$ with $h_p(a)= \min h_p$, thus $\delta(a)\in \BFra(s,r)$ under suitable congruence restrictions, either all $[q]$  in $\delta(a)$ are separated from each other by entries $[q+1]$ or vice versa. Hence these delta vectors are separable.

\begin{proposition}{\label{PropSeparable}}
Let $(r-1)\nmid (s-1)$, and let $g \equiv s-1 \bmod (r-1)$ be the least positive residue.
For any  $d\in \BFra(s,r)$ satisfying $h_p(\delta^{-1}(d))=\min h_p$, 
we have:
\begin{itemize}
\item[(i)] If $2g\geq r-1$, then $\eta_{\max}(d)=1$ and $q_1=1$, $q_2=\frac{g}{r-g-2}$.
\item[(ii)] If $2g \leq r-2$, then $\theta_{\max}(d)=1$ and $q_2=1$, $q_1=\frac{r-g-1}{g+1}$.
\end{itemize}
\end{proposition}

{\sc Proof. }\vspace{7pt}\newline
(i) Let $d=(d_1,\ldots,d_{r-1})$ and assume that $\eta_{\max}(d)\geq 2$. Hence there is some $1\leq u\leq r-2$ such that
$d_u=d_{u+1}=[q]$. Since $d_1=d_{r-1} =[q]$ and $g\geq r-g-1$, there is some $2\leq v \leq r-3$ such that $d_v=d_{v+1}= [q+1]$. By (\ref{symmhp}) we may assume that
$u<v$. Moreover, we may assume w.l.o.g. that $d_j\neq d_{j+1}$ for $u+1\leq j \leq v-1$ (otherwise we could choose $u$ larger or $v$ smaller, respectively). 
Thus, and since $d_j\in \{[q],[q+1]\}$, the sequence $(d_j)_{u+1\leq j \leq v}$  
is alternating, starting with $d_{u+1}=[q]$ and terminating with $d_v= [q+1]$.
Consequently, $v-u$ is even, and for $u+1\leq j \leq v$ we have
\[
d_j = \left\{ \begin{array}{cl}
                      \mbox{$[q]$}     & \mbox{ if $j\not\equiv u \bmod 2$, } \\
                      \mbox{$[q+1]$}  & \mbox{ if $j\equiv u \bmod 2$.}
                      \end{array} \right.
\]
This means that 
\begin{equation}\label{alternating}
a_{u+2+j} = a_{u+2} + [q]j +\left[\frac{j+1}{2}\right] \quad\quad (0\leq j \leq v-u-1).
\end{equation}

For $a=(a_1,\ldots,a_r):= \delta^{-1}(d)\in A(s,r)$, we define $b=(b_1,\ldots,b_r)\in A(s,r)$ by setting
\begin{equation}\label{vectorb}
 b_j := \left\{ \begin{array}{cl}
                  a_j & \mbox{ for $1\leq j\leq u+1$ or $v+1\leq j \leq r$, }\\
                  a_{j+1}-[q]& \mbox{ for $u+2\leq j \leq v$, }
                  \end{array} \right.
\end{equation}                        
i.e. we swap $d_{u+1}$ and $d_v$ in $d$.
Clearly, $\delta(b)\in \BFra(s,r)$.
Then, by (\ref{basishp}),
\begin{align}\label{diffhptotal}
\begin{split}
h_p(a)-h_p(b) 
&= \left(\sum_{k=1}^{u+1}\sum_{i=u+2}^v + \sum_{k=u+2}^{v-1}\sum_{i=k+1}^v + \sum_{k=u+2}^v\sum_{i=v+1}^r\right) \,\left(\frac{1}{p^{a_i-a_k}} - \frac{1}{p^{b_i-b_k}}\right)      \\
&= \sum_{k=1}^{u+1}\sum_{i=u+2}^v \,\left(\frac{1}{p^{a_i-a_k}} - \frac{1}{p^{(a_{i+1}-[q])-a_k}}\right)  \\
 &\quad\quad + \sum_{k=u+2}^{v-1}\sum_{i=k+1}^v \,\left(\frac{1}{p^{a_i-a_k}} - \frac{1}{p^{(a_{i+1}-[q])-(a_{k+1}-[q])}}\right) \\
 &\quad\quad + \sum_{k=u+2}^v\sum_{i=v+1}^r \,\left(\frac{1}{p^{a_i-a_k}} - \frac{1}{p^{a_i-(a_{k+1}-[q])}}\right)      \\
&= \sum_{k=1}^{u+1}\,p^{a_k}\sum_{i=u+2}^v \,\left(\frac{1}{p^{a_i}} - \frac{1}{p^{a_{i+1}-[q]}}\right)
 + \sum_{k=u+2}^{v-1}\sum_{i=k+1}^v \,\left(\frac{1}{p^{a_i-a_k}} - \frac{1}{p^{a_{i+1}-a_{k+1}}}\right) \\
 &\quad\quad + \sum_{k=u+2}^v\,\left(p^{a_k}-p^{a_{k+1}-[q]}\right) \sum_{i=v+1}^r \,\frac{1}{p^{a_i}}\, . 
\end{split}
\end{align}

For the middle double sum, we obtain
\begin{align}\label{hpmiddle1}
\begin{split}
\sum_{k=u+2}^{v-1}\sum_{i=k+1}^v \,\left(\frac{1}{p^{a_i-a_k}} - \frac{1}{p^{a_{i+1}-a_{k+1}}}\right) &= \sum_{k=u+2}^{v-1}  p^{a_k} \sum_{i=k+1}^v \,\frac{1}{p^{a_i}} - 
\sum_{k=u+2}^{v-1}  p^{a_{k+1}} \sum_{i=k+1}^v \,\frac{1}{p^{a_{i+1}}}  \\
&= \sum_{k=u+2}^{v-1}  p^{a_k} \sum_{i=k+1}^v \,\frac{1}{p^{a_i}} - 
\sum_{k=u+3}^v  p^{a_k} \sum_{i=k+1}^{v+1} \,\frac{1}{p^{a_i}}  \\
&= p^{a_{u+2}} \sum_{i=u+3}^v \,\frac{1}{p^{a_i}} - \frac{1}{p^{a_{v+1}}} \sum_{k=u+3}^v  p^{a_k}\\
&=  \sum_{j=1}^{v-u-2} \,\frac{1}{p^{a_{u+2+j}-a_{u+2}}} - \sum_{j=1}^{v-u-2} \,
\frac{1}{p^{a_{v+1}-a_{v+1-j}}}\,.
\end{split}
\end{align}
Now (\ref{alternating}) implies that 
$a_{u+2+j}-a_{u+2} = [q]j + \left[\frac{j+1}{2}\right] = a_{v+1}-a_{v+1-j}$
for $1\leq j \leq v-u-2$. Hence the last two sums in (\ref{hpmiddle1}) cancel termwise, and we conclude
\begin{equation}
\sum_{k=u+2}^{v-1}\sum_{i=k+1}^v \,\left(\frac{1}{p^{a_i-a_k}} - \frac{1}{p^{a_{i+1}-a_{k+1}}}\right) = 0. \label{hpmiddle2}
\end{equation}

By (\ref{alternating}), we also have
\begin{align}\label{part1}
\begin{split}
\sum_{i=u+2}^v \,\left(\frac{1}{p^{a_i}} - \frac{1}{p^{a_{i+1}-[q]}}\right) &=
\sum_{j=0}^{v-u-2} \,\left(\frac{1}{p^{a_{u+2+j}}} - \frac{1}{p^{a_{u+2+j+1}-[q]}}\right)     \\
&= \frac{1}{p^{a_{u+2}}} \sum_{j=0}^{v-u-2} \,\frac{1}{p^{[q]j}}\left(\frac{1}{p^{[\frac{j+1}{2}]}} -\frac{1}{p^{[\frac{j+2}{2}]}} \right)  \\
&= \frac{1}{p^{a_{u+2}}} \sum_{j=0}^{\frac{v-u}{2}-1}
\,\frac{1}{p^{(2[q]+1)j}}\left(1 -\frac{1}{p} \right)\\
&= (p-1)\frac{1}{p^{a_{u+2}+1}} \sum_{j=0}^{\frac{v-u}{2}-1}
\,\frac{1}{p^{(2[q]+1)j}}\,,
\end{split}
\end{align}
and similarly
\begin{align}\label{part2}
\begin{split}
\sum_{k=u+2}^v\,\left(p^{a_k}-p^{a_{k+1}-[q]}\right)  &=
p^{a_{u+2}} \sum_{j=0}^{v-u-2}\, p^{[q]j}\left(p^{[\frac{j+1}{2}]} -p^{[\frac{j+2}{2}]} \right) \\
&= p^{a_{u+2}} \sum_{j=0}^{\frac{v-u}{2}-1}\,p^{(2[q]+1)j}(1-p)\\
&= (1-p)p^{a_{u+2}+(2[q]+1)(\frac{v-u}{2}-1)} \sum_{j=0}^{\frac{v-u}{2}-1}\,\frac{1}{p^{(2[q]+1)j}}\,.
\end{split}
\end{align}
Using (\ref{hpmiddle2}), (\ref{part1}) and (\ref{part2}) altogether in (\ref{diffhptotal}) implies that
\begin{equation}\label{final1}
Q:=\frac{h_p(a)-h_p(b)}{(p-1)\sum\limits_{j=0}^{\frac{v-u}{2}-1}\,\frac{1}{p^{(2[q]+1)j}}}=  
\frac{1}{p^{a_{u+2}+1}} \sum_{k=1}^{u+1}\,p^{a_k} \; - \;\;\; 
 p^{a_{u+2}+(2[q]+1)(\frac{v-u}{2}-1)}\sum_{i=v+1}^r \,\frac{1}{p^{a_i}}\,. \\
\end{equation}
Since $a_{v+2}-a_{v+1} = d_{v+1} = [q+1]$ by definition, and since $d_j \geq [q]$ for all $j$, it follows that
\begin{align*}
\begin{split}
\sum_{i=v+1}^r \,\frac{1}{p^{a_i}} &= 
\frac{1}{p^{a_{v+1}}}\sum_{i=v+1}^r\, \frac{1}{p^{a_i-a_{v+1}}}  \\
&\leq \frac{1}{p^{a_{v+1}}} \left( 1 + \frac{1}{p^{[q+1]}} \sum_{i=0}^{\infty}
\frac{1}{p^{[q]i}} \right) =  \frac{1}{p^{a_{v+1}}} \left( 1 + \frac{1}{p^{[q+1]}-p} \right)\,.
\end{split}
\end{align*}
By (\ref{alternating}), we have $a_{v+1}-a_{u+2} = [q](v-u-1) + \frac{v-u}{2}$.
Applying this as well as the last inequality and $\sum_{k=1}^{u+1}\,p^{a_k} \geq p^{a_{u+1}} + p^{a_u}$ to (\ref{final1}), we obtain
\begin{align*}
Q &\geq p^{a_{u+1}-a_{u+2}-1} + p^{a_u -a_{u+2}-1} -p^{a_{u+2}-a_{v+1} +(2[q]+1)(\frac{v-u}{2}-1)}\left( 1 + \frac{1}{p^{[q+1]}-p} \right) \\
&= p^{-[q]-1} + p^{-2[q]-1} -p^{-[q]-1}\left( 1 + \frac{1}{p^{[q+1]}-p} \right) \\
&=  p^{-2[q]-1}\left( 1 - \frac{1}{p-p^{1-[q]}}\right)\,.
\end{align*}
This last term is positive because of $q\geq 1$. By definition of $Q$ in (\ref{final1}), we conclude that $h_p(a) > h_p(b)$. This contradicts the minimality condition for $h_p(a)$, and thus our initial assumption $\eta_{\max}(d)\geq 2$ must be wrong.
Therefore, $\eta_{\max}(d)=1$, which means that $q_1=\frac{r-g-1}{w+1}=1$. Hence $w=r-g-2$ and $q_2=\frac{g}{w}=\frac{g}{r-g-2}$.
\bigskip\newline
(ii) We assume that $\theta_{\max}(d)\geq 2$. Hence there is some $2\leq v\leq r-3$ such that
$d_v=d_{v+1}=[q+1]$. Since  $g\leq r-g-2$, there is some $1\leq u \leq r-2$ such that $d_u=d_{u+1}= [q]$. By (\ref{symmhp}) we may assume that
$u<v$. Moreover, we may assume w.l.o.g. that $d_j\neq d_{j+1}$ for $u+1\leq j \leq v-1$ (otherwise we could choose $u$ larger or $v$ smaller, respectively). At this point we are exactly in the same situation as in the proof of part (i). Again $b$ as defined in (\ref{vectorb}) reveals that $h_p(a)> \min h_p$, and this contradiction completes the proof of the proposition.

\qed

%
%
\section{Bivalence of second degree -- Proof of Theorem \ref{Thm2.2}}

We denote by $\BFSep(s,r)$ the set of all $d\in \BFra(s,r)$ having no neighbouring entries $[q]$ in case $2g\geq r-1$ and no neighbouring entries $[q+1]$ in case $2g \leq r-2$, respectively. Then Proposition \ref{PropSeparable} mainly says that $h_p(a)=\min h_p$ implies $\delta(a) \in \BFSep(s,r)$.

Assuming that $q\notin\mathbb{N}$,
we shall now see that in case $2g\ge r-1$ all $[q+1]$-blocks in $\delta(a)$, lying between two successive entries $[q]$, are of length $[q_2]$ or $[q_2+1]$. In case $2g\leq r-2$ all $[q]$-blocks in $\delta(a)$ have length either $[q_1]$ or $[q_1+1]$  (cf. (\ref{eta}) and (\ref{theta})).

If $(r-1)\nmid (s-1)$, we define for $d\in \BFSep(s,r)$
\[ \eta(d) := \eta_{\max}(d) - \eta_{\min}(d)\]
and 
\[ \theta(d) := \theta_{\max}(d) - \theta_{\min}(d).\]
i.e. $\eta(d)$ is the difference between the lengths of the longest and the shortest maximal sequence of successive values $[q]$ in $d$, and $\theta(d)$ is the corresponding difference for successive values $[q+1]$.
By Proposition \ref{PropSeparable} we know for any minimizer $a\in A(s,r)$ of $h_p$ that $\eta(\delta(a))=0$ in case $h\geq r-h-1$ and
$\theta(\delta(a))=0$ in case $h \leq r-h-2$. 

For a bivalent, separable integer vector $v$ we may formally derive a vector
$\Lambda(v)$ as follows. Let $m\neq k$ be the two entries of $v$ and assume w.l.o.g. 
that $v$ contains no consecutive entries $m$. If the same holds for $k$, then we assume $m<k$ for tie-breaking. Set $\Lambda(v):=(\lambda_1,\ldots,\lambda_{\ell})$ for suitable $\ell$,
where $\lambda_i$ is the length of the $i$-th maximal sequence of consecutive $k$-entries,
as separated by the $m$-entries. 
If $\Lambda(v)$, like $v$, is bivalent we shall call $v$ {\em bivalent of second degree}.

For $d\in \BFSep(s,r)$ we clearly have $\min\{\eta_{\max}(d),\theta_{\max}(d)\}=1$ due to
separability. The following proposition strengthens Proposition \ref{PropSeparable}
in the sense that, under the same assumptions on $r$ and $s$,
some $d\in\BFra(s,r)$ with $h_p(\delta^{-1}(d))= \min h_p$ is not only  separable
but also satisfies $\eta(d)+\theta(d) \le 1$. The latter amounts to the fact that $d$ is bivalent of second degree.

\begin{proposition}{\label{ProcBivSecond}}
Let $(r-1)\nmid (s-1)$, and let $g \equiv s-1 \bmod (r-1)$ be the least positive residue. If $d\in \BFra(s,r)$ satisfies $h_p(\delta^{-1}(d))= \min h_p$, then we have:
\begin{itemize}
\item[(i)] If $2g\geq r-1$, then $\eta_{\max}(d)=1$ and $\theta(d)\leq 1$.
\item[(ii)] If $2g \leq r-2$, then $\theta_{\max}(d)=1$ and $\eta(d)\leq 1$.
\end{itemize}
\end{proposition}

{\sc Proof. }\vspace{7pt} \newline
(i) Let $d=(d_1,\ldots,d_{r-1})$ satisfy the conditions of the proposition, in particular $d_1=d_{r-1}=[q]$, and there are integers 
$1=j_1 < j_2 < \ldots < j_{r-g-2} < j_{r-g-1}=r-1$  with the property
\[
d_j = \left\{ \begin{array}{cl}
                     \mbox{$[q]$}   & \mbox{ for $j\in \{j_1,j_2,\ldots,j_{r-g-1}\}$, } \\
                     \mbox{$[q+1]$} & \mbox{ for $j\notin \{j_1,j_2,\ldots,j_{r-g-1}\}$ } 
                     \end{array} \right.   \quad\quad\quad (1\leq j \leq r-1).
\]
It follows from Proposition \ref{PropSeparable}(i) that $\eta_{\max}(d)=1$, hence
$j_{i+1}-j_i \geq 2$ for $1\leq i \leq r-g-2$. 

In order to prove the other assertion of (i) we make the assumption that $\theta(d)\geq 2$, i.e. there are two groups of successive entries $[q+1]$ in $d$ whose lengths differ by at least $2$. 
Hence, using the notation introduced in (\ref{deftj}), we can find integers $1\leq u\leq w$ and $1\leq v \leq w$ such that $j_{u+1} - j_u-1 = t_{2u}$ and $j_{v+1} - j_v-1 = t_{2v}$ satisfy $t_{2v}-t_{2u} \geq 2$, and we may assume that $|v-u|$ is minimal with this property. By (\ref{symmhp}) we can also assume w.l.o.g. that $u<v$. We therefore have
\begin{align*}
d &=     
(\ldots,d_{j_u},\underbrace{[q+1],\ldots,[q+1]}_{t_{2u}-fold},d_{j_{u+1}},\ldots\ldots,d_{j_v},\underbrace{[q+1],\ldots,[q+1]}_{t_{2v}-fold},d_{j_{v+1}}\ldots)\\
&= (\ldots,[q],\underbrace{[q+1],\ldots,[q+1]}_{t_{2u}-fold},[q],\ldots\ldots,[q], \underbrace{[q+1],\ldots,[q+1]}_{t_{2v}-fold},[q]\ldots),  
\end{align*}
and the desired contradiction will be derived in two steps: We first deal with the case where merely a single $[q]$-block separates the two $[q+1]$-blocks of lengths $t_{2u}$ and $t_{2v}$, and later we shall handle greater distances between them. In both situations, we construct some $b\in A(s,r)$ satisfying $\delta(b) \in \BFra(s,r)$ and $h_p(b)<h_p(a)$ by counterbalancing the lengths of the two $[q+1]$-blocks. We set 
$a=(a_1,\ldots,a_r):= \delta^{-1}(d)$.

\underline{Case 1}: $v=u+1$.
\newline
It follows that
\begin{equation}\label{ajot1}
a_{j_v+1} -a_k = \left\{
                 \begin{array}{ll}
                 (j_v-k+1)[q+1]-1 & \mbox{ for $j_u+1\leq k \leq j_v$, } \\
                 (j_v-j_u+1)[q+1]-2  & \mbox{ for $k=j_u$, }
                 \end{array}
                 \right.
\end{equation}
and
\begin{equation}\label{ajot2}
a_i - a_{j_v+1} =  (i-j_v-1)[q+1] \quad\quad (j_v+2 \leq i \leq j_{v+1})\,. 
\end{equation}

We define $b=(b_1,\ldots,b_r)\in A(s,r)$ by setting
\begin{equation}\label{vectorb1}
 b_j := \left\{ \begin{array}{cl}
                  a_j & \mbox{ for $1\leq j\leq j_v$ or $j_v+2\leq j \leq r$, }\\
                  a_{j_v+1}+1 & \mbox{ for $j=j_v+1$. }
                  \end{array} \right.
\end{equation}   
Clearly, $\delta(a)=d \in \BFra(s,r)$ implies $\delta(b) \in \BFra(s,r)$.  
Then, by (\ref{basishp}),
\begin{align*}
\begin{split}
h_p(a)-h_p(b) 
&= \sum_{i=j_v+2}^r \,\left(\frac{1}{p^{a_i-a_{j_v+1}}} - \frac{1}{p^{b_i-b_{j_v+1}}}\right)
 + \sum_{k=1}^{j_v}\,\left(\frac{1}{p^{a_{j_v+1}-a_k}} - \frac{1}{p^{b_{j_v+1}-b_k}}\right)      \\
&=  \sum_{i=j_v+2}^r \, \frac{1}{p^{a_i}}\left(p^{a_{j_v+1}} - p^{a_{j_v+1}+1}\right) 
+ \sum_{k=1}^{j_v}\,p^{a_k} \left(\frac{1}{p^{a_{j_v+1}}} - \frac{1}{p^{a_{j_v+1}+1}}\right)      \\
&= (p-1)\left(\sum_{k=1}^{j_v}\,\frac{1}{p^{a_{j_v+1}-a_k+1}} -   \sum_{i=j_v+2}^r\, \frac{1}{p^{a_i-a_{j_v+1}}}\right) . 
\end{split}
\end{align*}
We obtain
\begin{align*}
\begin{split}
\frac{h_p(a)-h_p(b)}{p-1} 
&> \frac{1}{p^{a_{j_v+1}-a_{j_u}+1}}+ \sum_{k=j_u+1}^{j_v}\,\frac{1}{p^{a_{j_v+1}-a_k+1}} \\
&\quad\quad -   \sum_{i=j_v+2}^{2j_v-j_u+1}\, \frac{1}{p^{a_i-a_{j_v+1}}} -
\sum_{i=2j_v-j_u+2}^{\infty}\, \frac{1}{p^{a_i-a_{j_v+1}}} 
\end{split}
\end{align*}
and observe that 
the first two sums on the righthand side have the same number of terms. Since
\begin{equation} 
2j_v-j_u+2 = j_v+j_{u+1}-j_u+2 = j_v + t_{2u} +3 \leq j_v + t_{2v} +1 =j_{v+1},
\label{ajot3}
\end{equation}
we can apply (\ref{ajot1}) and (\ref{ajot2}) to deduce termwise cancellation of those two sums. Hence, and by (\ref{ajot1}), (\ref{ajot3}) and (\ref{ajot2}) again, it follows that
\begin{align*}
\begin{split}
\frac{h_p(a)-h_p(b)}{p-1} 
&> \frac{1}{p^{a_{j_v+1}-a_{j_u}+1}} -
\sum_{i=2j_v-j_u+2}^{\infty}\, \frac{1}{p^{a_i-a_{j_v+1}}} \\
&\geq \frac{1}{p^{a_{j_v+1}-a_{j_u}+1}} -
  \frac{1}{p^{a_{2j_v-j_u+2}-a_{j_v+1}}}     \sum_{i=0}^{\infty}\, \frac{1}{p^i} \\
&= \frac{1}{p^{(j_v-j_u+1)[q+1]-1}} - \frac{1}{p^{(j_v-j_u+1)[q+1]}}\, \frac{p}{p-1} \\
&= \frac{1}{p^{(j_v-j_u+1)[q+1]-1}}\left( 1 - \frac{1}{p-1}\right) \geq 0\,.
\end{split}
\end{align*}
This contradicts the minimality condition for $h_p(a)$, and thus our initial assumption $\theta(d)\geq 2$ must be wrong in this case.

\underline{Case 2}: $v\geq u+2$.
\newline
By Case 1, we know that 
\[ t_{2v} - t_{2u} = ( t_{2v} - t_{2(u+1)}) + (t_{2(u+1)} - t_{2u}) \leq 
\vert t_{2v} - t_{2(u+1)}\vert +1.\]
Now the assumption $t_{2v}-t_{2u} \geq 3$ would imply $\vert t_{2v} - t_{2(u+1)}\vert \geq 2$, contradicting the minimality of $\vert v-u\vert$. 
We are left with $t_{2v}-t_{2u} = 2$. The minimality of $\vert v-u\vert$ implies in this special situation that 
\[ 
t_{2u}+1 =t_{2(u+1)} = t_{2(u+2)} = \ldots = t_{2(v-1)} =  t_{2v}-1,
\]
i.e. we have for $\Delta_u(d):= j_{u+1} -j_u+1$ that
\begin{equation}
\Delta_u(d)= j_{u+2}-j_{u+1} = j_{u+3}-j_{u+2} = \ldots = j_v-j_{v-1}= j_{v+1}-j_v-1.
\label{grossdelta}
\end{equation}
We also have
\begin{equation}\label{ajot11}
a_{j_{u+1}+1} -a_k = \left\{
                 \begin{array}{ll}
                 (j_{u+1}-k+1)[q+1]-1 & \mbox{ for $j_u+1\leq k \leq j_{u+1}$, } \\
                 (j_{u+1}-j_u+1)[q+1]-2  & \mbox{ for $k=j_u$, }
                 \end{array}
                 \right.
\end{equation}
and
\begin{equation}\label{ajot12}
a_i - a_{j_v+1} =  (i-j_v-1)[q+1] \quad\quad (j_v+2 \leq i \leq j_{v+1})\,. 
\end{equation}

We define $b=(b_1,\ldots,b_r)\in A(s,r)$  by setting
\begin{equation}\label{vectorb11}
 b_j := \left\{ \begin{array}{cl}
                  a_j & \mbox{ for $1\leq j\leq j_{u+1}$ or $j_v+2\leq j \leq r$, }\\
                  a_{j-1}+[q+1]& \mbox{ for $j_{u+1}+1\leq j \leq j_v+1$. }
                  \end{array} \right.
\end{equation}                        
i.e. we enlarge the number of intervals of length $[q+1]$ between $a_{j_u+1}$ and $a_{j_{u+1}}$ by one and shorten the number of these intervals between $a_{j_v+1}$ and $a_{j_{v+1}}$ by one. Clearly, $\delta(a)=d \in \BFra(s,r)$ implies $\delta(b) \in \BFra(s,r)$.
Then, by (\ref{basishp}),
\begin{align}\label{diffhp11}
\begin{split}
h_p(a)-h_p(b) 
&= \left(\sum_{k=1}^{j_{u+1}}\sum_{i=j_{u+1}+1}^{j_v+1} + \sum_{k=j_{u+1}+1}^{j_v+1}\sum_{i=k+1}^{j_v+1} + \sum_{k=j_{u+1}+1}^{j_v+1}\sum_{i=j_v+2}^r\right) \,\left(\frac{1}{p^{a_i-a_k}} - \frac{1}{p^{b_i-b_k}}\right)      \\
&= \sum_{k=1}^{j_{u+1}}\sum_{i=j_{u+1}+1}^{j_v+1} \,\left(\frac{1}{p^{a_i-a_k}} - \frac{1}{p^{(a_{i-1}+[q+1])-a_k}}\right)  \\
 &\quad\quad + \sum_{k=j_{u+1}+1}^{j_v}\sum_{i=k+1}^{j_v+1} \,\left(\frac{1}{p^{a_i-a_k}} - \frac{1}{p^{(a_{i-1}+[q+1])-(a_{k-1}+[q+1])}}\right) \\
 &\quad\quad + \sum_{k=j_{u+1}+1}^{j_v+1}\sum_{i=j_v+2}^r \,\left(\frac{1}{p^{a_i-a_k}} - \frac{1}{p^{a_i-(a_{k-1}+[q+1])}}\right)      \\
&= \sum_{k=1}^{j_{u+1}}\, p^{a_k} \sum_{i=j_{u+1}+1}^{j_v+1}  
\,\left(\frac{1}{p^{a_i}} - \frac{1}{p^{a_{i-1}+[q+1]}}\right) \\
&\quad\quad +  \sum_{k=j_{u+1}+1}^{j_v}\sum_{i=k+1}^{j_v+1}\,\left(\frac{1}{p^{a_i-a_k}} - \frac{1}{p^{a_{i-1}-a_{k-1}}}\right) \\
 &\quad\quad +\sum_{k=j_{u+1}+1}^{j_v+1} \,\left(p^{a_k}-p^{a_{k-1}+[q+1]}\right) \sum_{i=j_v+2}^r\,\frac{1}{p^{a_i}}\, . 
\end{split}
\end{align}
By definition of the $j_{\ell}$, we have for $j_{u+1}+1 \leq i \leq j_v+2-1$
\begin{equation}
a_i -a_{i-1} = d_{i-1} = \left\{ \begin{array}{cl}
                  $[$q$]$ & \mbox{ for $i=j_{\ell}+1$, $u+1\leq \ell \leq v$, }\\
                  $[$q+1$]$ & \mbox{ otherwise. }
                  \end{array} \right.
\label{diffai}
\end{equation}
This implies that
\begin{equation}
\sum_{i=j_{u+1}+1}^{j_v+1}  \,\left(\frac{1}{p^{a_i}} - \frac{1}{p^{a_{i-1}+[q+1]}}\right)
= \sum_{\ell=u+1}^v \,\left(\frac{1}{p^{a_{j_{\ell}+1}}} - \frac{1}{p^{a_{j_{\ell}+1}+1}}\right)
= \left(1-\frac{1}{p}\right) \sum_{\ell=u+1}^v \,\frac{1}{p^{a_{j_{\ell}+1}}}
\label{inter1}
\end{equation}
and
\begin{equation}
\sum_{k=j_{u+1}+1}^{j_v+1}  \,\left(p^{a_k}-p^{a_{k-1}+[q+1]}\right)
= \sum_{\ell=u+1}^v \,\left(p^{a_{j_{\ell}+1}} - p^{a_{j_{\ell}+1}+1}\right)
= \left(1-p\right) \sum_{\ell=u+1}^v \,p^{a_{j_{\ell}+1}}.
\label{inter2}
\end{equation} 
Moreover
\begin{align}\label{diffmiddle}
\begin{split}
\sum_{k=j_{u+1}+1}^{j_v}\sum_{i=k+1}^{j_v+1}\,&\left(\frac{1}{p^{a_i-a_k}} - \frac{1}{p^{a_{i-1}-a_{k-1}}}\right) \\
&= \sum_{k=j_{u+1}+1}^{j_v} p^{a_k} \sum_{i=k+1}^{j_v+1}\,\frac{1}{p^{a_i}} -
  \sum_{k=j_{u+1}}^{j_v-1} p^{a_k} \sum_{i=k+1}^{j_v}\,\frac{1}{p^{a_i}}  \\
&=  \frac{1}{p^{a_{j_v+1}}} \sum_{k=j_{u+1}+1}^{j_v} p^{a_k} - 
   p^{a_{j_{u+1}}} \sum_{i=j_{u+1}+1}^{j_v}\, \frac{1}{p^{a_i}}\, \\
&= \sum_{k=j_{u+1}+1}^{j_v} \frac{1}{p^{a_{j_v+1}-a_k}} - \sum_{i=j_{u+1}+1}^{j_v}\, \frac{1}{p^{a_i-a_{j_{u+1}}}} \\
&= \sum_{k=1}^{j_v-j_{u+1}} \frac{1}{p^{a_{j_v+1}-a_{j_v+1-k}}} - \sum_{i=1}^{j_v-j_{u+1}}\, \frac{1}{p^{a_{j_{u+1}+i}-a_{j_{u+1}}}}\, .
\end{split}
\end{align}
It is easy to deduce from (\ref{diffai}) and (\ref{grossdelta}) that, by symmetry of the spacing,  
\[ a_{j_v+1} -a_{j_v+1-k}  
= a_{j_{u+1}+k} -a_{j_{u+1}} \]
for $0 \leq k \leq j_v- j_{u+1} +1$. Hence the last two sums in (\ref{diffmiddle}) cancel termwise, and we obtain from (\ref{diffhp11}), (\ref{inter1}) and (\ref{inter2}) that
\begin{align}\label{diffhp15}
\begin{split}
\frac{h_p(a)-h_p(b)}{p-1} &= \frac{1}{p} \sum_{k=1}^{j_{u+1}}\, p^{a_k} \sum_{\ell=u+1}^v \,\frac{1}{p^{a_{j_{\ell}+1}}} - \sum_{i=j_v+2}^r\,\frac{1}{p^{a_i}}
 \sum_{\ell=u+1}^v \,p^{a_{j_{\ell}+1}}  \\
 &= \frac{1}{p^{a_{j_{u+1}+1}+1}} \sum_{k=1}^{j_{u+1}}\, p^{a_k} \sum_{\ell=u+1}^v \,\frac{1}{p^{a_{j_{\ell}+1}-a_{j_{u+1}+1}}}   \\
 &\quad\quad\quad  - p^{a_{j_v+1}} \sum_{i=j_v+2}^r\,\frac{1}{p^{a_i}}
 \sum_{\ell=u+1}^v \,\frac{1}{p^{a_{j_v+1}-a_{j_{\ell}+1}}} \\
&= \frac{1}{p^{a_{j_{u+1}+1}+1}} \sum_{k=1}^{j_{u+1}}\, p^{a_k} \sum_{\ell=0}^{v-u-1} \,\frac{1}{p^{(\Delta_u(d)[q+1]-1)\ell}} \\
 &\quad\quad\quad   - p^{a_{j_v+1}} \sum_{i=j_v+2}^r\,\frac{1}{p^{a_i}}
 \sum_{\ell=0}^{v-u-1} \,\frac{1}{p^{(\Delta_u(d)[q+1]-1)\ell}} \\ 
 &= \sum_{\ell=0}^{v-u-1} \,\frac{1}{p^{(\Delta_u(d)[q+1]-1)\ell}}  \left( \sum_{k=1}^{j_{u+1}}\, \frac{1}{p^{a_{j_{u+1}+1}-a_k+1}}  -  \sum_{i=j_v+2}^r\,\frac{1}{p^{a_i-a_{j_v+1}}} \right). 
\end{split}
\end{align}
 
We have
\begin{align*}
\begin{split}
\sum_{k=1}^{j_{u+1}}\, &\frac{1}{p^{a_{j_{u+1}+1}-a_k+1}}  -  
 \sum_{i=j_v+2}^r\,\frac{1}{p^{a_i-a_{j_v+1}}} \\
&\quad\quad\quad>  \frac{1}{p^{a_{j_{u+1}+1}-a_{j_u}+1}}+ \sum_{k=j_u+1}^{j_{u+1}}\,\frac{1}{p^{a_{j_{u+1}+1}-a_k+1}} \\
&\quad\quad\quad\quad  -   \sum_{i=j_v+2}^{j_v+j_{u+1}-j_u+1}\, \frac{1}{p^{a_i-a_{j_v+1}}} -
\sum_{i=j_v+j_{u+1}-j_u+2}^{\infty}\, \frac{1}{p^{a_i-a_{j_v+1}}} 
\end{split}
\end{align*}
and observe that the positive sum and the first negative sum on the righthand side have the same number of terms. Since, by (\ref{grossdelta}), 
\begin{equation} 
j_v+j_{u+1}-j_u+2 = j_{v+1}-(j_{v+1}-j_v)+(j_{u+1}-j_u)+2  =j_{v+1},
\label{ajot13}
\end{equation}
we can apply (\ref{ajot11}) and (\ref{ajot12}) to deduce termwise cancellation of those two sums. Hence, and by (\ref{ajot11}), (\ref{ajot13}) and (\ref{ajot12}) again, it follows that
\begin{align*}
\begin{split}
\sum_{k=1}^{j_{u+1}}\, \frac{1}{p^{a_{j_{u+1}+1}-a_k+1}}  -  
 \sum_{i=j_v+2}^r\,\frac{1}{p^{a_i-a_{j_v+1}}}
&>  \frac{1}{p^{a_{j_{u+1}+1}-a_{j_u}+1}}
- \sum_{i=j_v+j_{u+1}-j_u+2}^{\infty}\, \frac{1}{p^{a_i-a_{j_v+1}}} \\
&\geq  \frac{1}{p^{a_{j_{u+1}+1}-a_{j_u}+1}}  -
  \frac{1}{p^{a_{j_{v+1}}-a_{j_v+1}}}     \sum_{i=0}^{\infty}\, \frac{1}{p^i} \\
&=  \frac{1}{p^{(j_{u+1}-j_u+1)[q+1]-1}}  -
  \frac{1}{p^{(j_{v+1}-j_v-1)[q+1]}}  \, \frac{p}{p-1}    \\ 
&=  \frac{1}{p^{\Delta_u(d)[q+1]-1}}  -
  \frac{1}{p^{\Delta_u(d)[q+1]}}  \, \frac{p}{p-1}    \\ 
&= \frac{1}{p^{\Delta_u(d)[q+1]-1}}\left( 1 - \frac{1}{p-1}\right) \geq 0\,.
\end{split}
\end{align*}
With this inequality, (\ref{diffhp15}) implies $h_p(a)>h_p(b)$, contradicting the minimality condition for $h_p(a)$. Again the initial assumption $\theta(d) \geq 2$ cannot hold, which completes the proof of (i).
\bigskip\newline
(ii) For $2g \leq r-2$, it follows from Proposition \ref{PropSeparable}(ii) that $\theta_{\max}(d)=1$. By the respective arguments, corresponding directly to the ones used in (i), now the assumption $\eta(d)\geq 2$ turns out to be contradictive.

\qed

\bigskip
{\sc Proof of Theorem \ref{Thm2.2}. }\newline 
(i) We know from Proposition \ref{PropFraming} that $d=(d_1,\ldots,d_{r-1}):=\delta(a)\in \BFra(s,r)$.
Since $(r-g-2)\mid g$, the number $q_2$ is an integer. It follows from (\ref{sumstl}) that the existence of a $t_{2k}< q_2$ would imply the existence of a $t_{2\ell}> q_2$ and vice versa, both cases contradicting $\theta_{\max}(d)-\theta_{\min}(d)=\theta(d)\leq 1$, which holds by Proposition \ref{ProcBivSecond}(i). Hence  $t_{2\ell}= q_2$ for all $\ell$, and the assertion follows. 
\bigskip\newline
(ii) The argument in the proof of (i) showed that the existence of a $t_{2k}< q_2$ implies the existence of a $t_{2\ell} > q_2$ and vice versa. Hence $\theta_{\min}(d)=[q_2]$ and $\theta_{\max}(d)=[q_2]+1$. Denote by $x$ the number of sequences of $[q+1]$-blocks of length $[q_2]$ in $d$. Then there are $r-g-2-x$ sequences of $[q+1]$-blocks of length $[q_2]+1$ in $d$. It follows that
\[ x[q_2] + (r-g-2-x)([q_2]+1) =g,\]
hence $x=(r-g-2)[q_2+1]-g$. By definition, $(r-g-2)q_2 = g$ and $e= (r-g-2)(q_2-[q_2])$. These identities imply
\[ x -(r-g-2) = (r-g-2)[q_2] -g =(r-g-2)q_2 -e -g = -e,\]
which completes the proof of (ii).
\bigskip\newline
(iii) can be shown by the same reasoning as (i).
\bigskip\newline
(iv) follows like (ii).

\qed

%
%

\section{Continued balancing}

Recall that we cited an analytical result from \cite{SA2}
stating that $h_p(a)$ becomes minimal if $0=a_1<a_2<\ldots<a_{r-1}<a_r=s-1$ are chosen in nearly equidistant position. In terms of delta vector structure, we now see that bivalence is the first balancing step towards
this goal. Further balancing is achieved by placing the rarer of the two elements of the delta vector 
as singletons. This is the separability property. Finally, we expect that the separating 
singletons are again distributed in nearly equidistant position, which amounts to
bivalence of second degree. 

\medskip

The following example demonstrates the balancing effect numerically.

\begin{example} Let $s=22$ and $r=17$. Then 
\[\delta(a_1)=(1, 1, 2, 1, 1, 2, 1, 1, 2, 1, 2, 1, 1, 2, 1, 1)\]
gives the minimal possible value of $h_3(a_1)\approx 5.36266$, thus maximizing the energy among
all tuples of $A(s,r)$. On the other hand, the vector
\[\delta(a_2)=(1, 1, 1, 1, 1, 1, 1, 1, 1, 1, 1, 1, 1, 1, 1, 6)\]
gives a particularly large value of $h_3(a_2)\approx 7.25206$.

Restricting ourselves to bivalent delta vectors, the maximal value of $h_3$ achievable
is $h_3(a_3)\approx  5.96811$ for
\[\delta(a_3)=(2, 2, 2, 1, 1, 1, 1, 1, 1, 1, 1, 1, 1, 1, 2, 2).\]

A further restriction to delta vectors that are also $[q]$-framed yields
a maximal $h_3(a_4) \approx 5.79688$ for
\[\delta(a_4)=(1, 2, 2, 2, 2, 2, 1, 1, 1, 1, 1, 1, 1, 1, 1, 1).\]

If we additionally impose separability we get a maximal $h_3(a_5) \approx 5.47795$ for
\[\delta(a_5)=(1, 2, 1, 2, 1, 2, 1, 1, 1, 1, 1, 1, 2, 1, 2, 1).\]

Finally also requiring bivalence of second degree, we arrive at a maximal
$h_3(a_6) \approx 5.37484$ for
\[\delta(a_6)=(1, 2, 1, 1, 2, 1, 1, 2, 1, 1, 2, 1, 1, 2, 1, 1).\]

We can see that this is now quite close to $h_3(a_1) \approx 5.36266$.
\end{example}

In view of this example one tends to expect that the balancing continues as far as possible,
finally resulting in the desired energy maximizing divisor set.

A definition of balancing of a certain degree is readily derived.
Let us formally define 
\[
\begin{split}
\Lambda^0(d) &=d,\\
\Lambda^{i}(d)&=\Lambda(\Lambda^{i-1}(d)), \quad \mbox{for}~i\in\NN
\end{split}
\]
and say that $d$ 
is {\em balanced of $i$-th degree} if
$\Lambda^{i}(d)$ exists, i.e.~$\Lambda^{i-1}(d)$ is bivalent and separable.
We exclude cases where $\Lambda^{i}(d)$ would formally exist but be an empty vector due to
$\Lambda^{i-1}(d)$ having only identical entries.
Let us call $\Lambda^{0}(d),\ldots,\Lambda^{j}(d)$ the {\em $\Lambda$ sequence of $d$}
if $\Lambda^{i}(d)$ exists for $i=0,\ldots,j$, but not for $i=j+1$.
\medskip

The effect of continued balancing becomes strikingly apparent in the following example, where 
we have several levels of balancing.

\begin{example}\label{s44r35ex}
For $s=44$ and $r=35$ the $\Lambda$ sequence of the energy maximal divisor tuple (up to symmetry of the delta vector) is
\[
\begin{split}
&(1, 1, 1, 2, 1, 1, 2, 1, 1, 1, 2, 1, 1, 2, 1, 1, 2, 1, 1, 1, 2, 1, 1, 2, 1, 1, 1, 2, 1, 1, 2, 1, 1, 1), \\
&(3, 2, 3, 2, 2, 3, 2, 3, 2, 3),\\
&(1, 2, 1, 1),\\
&(1, 2), \\
&(1).
\end{split}
\]
\end{example}

It seems that framing is an important aspect in continued balancing.

\begin{conjecture}
Let $a$ be an energy maximal exponent tuple. Suppose that for $d:=\delta(a)$ all 
$\Lambda^0(d),\ldots,\Lambda^j(d)$ exist and that $\Lambda^j(d)$ is unframed.
Then $\Lambda^0(d),\ldots,\Lambda^{j-1}(d)$ are all framed.
\end{conjecture}

Formally, a vector $\Lambda^i(d)$ can be interpreted as the delta vector
$\delta(a')$ of some unique admissible exponent tuple $a'$ (with length $r'$ and largest entry $s'-1$). 
In this sense, one could shorten a given sequence $\Lambda^0(a),\ldots,\Lambda^j(a)$
to obtain the tail $\Lambda^0(a'),\ldots,\Lambda^{j-i}(a')$. This is shown in the next example.

\begin{example}\label{s26r11ex}
Based on the $\Lambda$ sequence given in Example \ref{s44r35ex}, consider the following vectors:

\[
\begin{split}
d=&(1, 1, 1, 2, 1, 1, 2, 1, 1, 1, 2, 1, 1, 2, 1, 1, 2, 1, 1, 1, 2, 1, 1, 2, 1, 1, 1, 2, 1, 1, 2, 1, 1, 1), \\
a=&(0,1,2,3,5,6,7,9,10,11,\ldots,32,33,34,36,37,38,40,41,42,43),\\
d'=&(3, 2, 3, 2, 2, 3, 2, 3, 2, 3),\\
a'=&(0,3,5,8,10,12,15,17,20,22,25).
\end{split}
\]

\medskip

Shortening the $\Lambda$ sequence of $d$ with $a=\delta^{-1}(d)$ by omitting the first delta vector
gives the $\Lambda$ sequence of $d'=\Lambda(d)$ with admissible exponent tuple $a'=\delta^{-1}(d')$, in which case we have $s'=26$ and $r'=11$:
\[
\begin{split}
&(3, 2, 3, 2, 2, 3, 2, 3, 2, 3),\\
&(1, 2, 1, 1),\\
&(1, 2), \\
&(1).
\end{split}
\]
\end{example}

It would be a most desirable property if $a$ were energy maximal within $A(r,s)$ that
the same would hold for $a'$ within $A(r',s')$. Examples indicate that this is often the case,
but not in general. Consider the following example:

\begin{example}
Consider the $\Lambda$ sequence given in Example \ref{s26r11ex}.

Clearly, the vector $(3, 2, 3, 2, 2, 3, 2, 3, 2, 3)$ does not define an energy maximal
divisor tuple since it does not have the $[q]$-framing property required by  Proposition \ref{PropFraming}.

\medskip

And indeed, the $\Lambda$ sequence of the energy maximal divisor tuple (again, up to symmetry) is
\[
\begin{split}
&(2, 3, 2, 3, 3, 2, 3, 2, 3, 2),\\
&(1, 2, 1, 1),\\
&(1, 2),\\
&(1).
\end{split}
\]
\end{example}

Although a continued balancing with longest possible sequences $\Lambda^0(d),\ldots,\Lambda^j(d)$
yields divisor tuples $a:=\delta^{-1}(d)\in A(s,r)$ with high energy, it does not automatically guarantee maximal energy
among the elements of $A(s,r)$. This can be seen from the next example. However, we suspect that
this effect is due to a probably not yet completely suitable formal notion of continued balancing. 

\begin{example}
For $s=16$ and $r=12$ the $\Lambda$ sequence of the energy maximal divisor tuple (up to symmetry of the delta vector) is
\[
\begin{split}
&(1, 1, 2, 1, 2, 1, 2, 1, 2, 1, 1), (2, 1, 1, 1, 2), (3)
\end{split}
\]
but the $\Lambda$ sequence of the runner-up is longer:
\[
\begin{split}
& (1, 2, 1, 1, 2, 1, 2, 1, 2, 1, 1), (1, 2, 1, 1, 2),  (1, 2), (1).
\end{split}
\]

Interestingly, this situation is reversed for $s=16$ and $r=11$:
\[
\begin{split}
&(1, 2, 1, 2, 1, 2, 2, 1, 2, 1), (1, 1, 2, 1), (2, 1), (1)
\end{split}
\]
is the $\Lambda$ sequence of the energy maximal divisor tuple, whereas
\[
\begin{split}
&(1, 2, 1, 2, 1, 2, 1, 2, 2, 1), (1, 1, 1, 2), (3)
\end{split}
\]
is the $\Lambda$ sequence of the runner-up.

\medskip

Note that in the first case we have $2g \le r-2$ and in the second case
$2g \ge r-1$. So, in view of the cases listed in Theorem \ref{Thm2.2},
we have a notable difference here that may have to do with the effect.
\end{example}

To better understand this process and properly embed it in a theory would be the object of
future work. In this context, let us remark that the continued balancing somewhat resembles what
happens in leap year calculations, which in turn are related to the Bresenham line drawing algorithm,
continued fractions and the Euclidean algorithm (cf. \cite{HAR}). Balancing also seems to be reminiscent of Beatty sequences and the way they partition $\ZZ$ into two sets (cf. \cite{ST}). Successfully linking these concepts
with maximizing the energy of integral ciculant graphs of prime power order is certainly a goal
inviting further research.

%
%


\begin{thebibliography}{999}
\addcontentsline{toc}{}{References}


\bibitem{AHM}
O. Ahmadi and N. Alon and I.F. Blake and I.E. Shparlinski,
Graphs with integral spectrum,
Linear Algebra Appl. {\bf 430} (2009), 547-552.

\bibitem{AKH} R. Akhtar and M. Boggess and T. Jackson-Henderson and I. Jim{\'e}nez and R. Karpman and A. Kinzel and D. Pritikin,
On the unitary Cayley graph of a finite ring,   
Electron. J. Combin. {\bf 16} (2009), Research Paper R117, 13 pp. (electronic). 

\bibitem{BAS} M. Ba{\v s}i{\'c} and A. Ili{\'c},
{On the clique number of integral circulant graphs},
Appl. Math. Lett. {\bf 22} (2009), 1406-1411.

\bibitem{BAS2} M. Ba{\v s}i{\'c} and A. Ili{\'c},
{On the Automorphism Group of Integral Circulant Graphs},
Electron. J. Combin. {\bf 18} (2011), Research Paper P68, 13 pp. (electronic). 

\bibitem{BAS3} M. Ba\v{s}i\'c and M.D. Petkovi\'c,
{Perfect state transfer in integral circulant graphs of non-square-free order},
Linear Algebra Appl. {\bf 433} (2010), 149-163

\bibitem{BAS4} M. Ba\v si\'c and M.D. Petkovi\'c and D. Stevanovi\'c,
{Perfect state transfer in integral circulant graphs},
Appl. Math. Lett. {\bf 22} (2009), 1117-1121

\bibitem{BEA} N. de Beaudrap,
{On restricted unitary {C}ayley graphs and symplectic transformations modulo {$n$}},
Electron. J. Combin. {\bf 17} (2010), Research Paper R69, 26 pp. (electronic). 

\bibitem{BER} P. Berrizbeitia and R. E. Giudici,
{On cycles in the sequence of unitary Cayley graphs},
Discrete Math. {\bf 282} (2004), 239-243.

\bibitem{BRU} R.A. Brualdi,
Energy of a graph, AIM Workshop Notes, 2006.

\bibitem{DAV} P. J. Davis,
{Circulant matrices},
John Wiley \&\ Sons, New York-Chichester-Brisbane, 1979.
  
\bibitem{DEJ} I. J. Dejter and R. E. Giudici,
{On unitary Cayley graphs},
J. Combin. Math. Combin. Comput. {\bf 18} (1995), 121-124.
   
\bibitem{DRO} A. Droll,
{A classification of {R}amanujan unitary {C}ayley graphs},
Electron. J. Combin. {\bf 17} (2010), Research Note N29, 6 pp. (electronic). 

\bibitem{GUT} I. Gutman,
{The energy of a graph},
Ber. Math.-Stat. Sekt. Forschungszent. Graz {\bf 103}, 1978.

\bibitem{HAR} M.A. Harris and E.M. Reingold,
Line {D}rawing, {L}eap {Y}ears, and {E}uclid,
ACM Computing Surveys {\bf 36} (2004), 68-80.
   
\bibitem{ILI} A. Ili\'{c},
The energy of unitary Cayley graphs,
Lin. Alg. Appl. {\bf 431} (2009), 1881-1889.

\bibitem{KLO} W. Klotz and T. Sander,
Some properties of unitary Cayley graphs,
Electron. J. Combin. {\bf 14} (2007), Research Paper R45, 12 pp. (electronic).

\bibitem{KLO2} W. Klotz and T. Sander,
{Integral {C}ayley graphs over abelian groups},
Electron. J. Combin. {\bf 17} (2010), Research Paper R81, 13 pp. (electronic). 

\bibitem{PET} M.D. Petkovi\'c and M. Ba\v{s}i\'c,
{Further results on the perfect state transfer in integral circulant graphs},
Comput. Math. Appl. {\bf 61} (2011), 300-312

\bibitem{RAM} H.N. Ramaswamy and C.R. Veena,
{On the Energy of Unitary Cayley Graphs},
Electron. J. Combin. {\bf 16} (2009), Research Note N24, 8 pp. (electronic).

\bibitem{SA1} J.W. Sander and T. Sander,
{The energy of integral circulant graphs with prime power order},
Appl. Anal. Discrete Math. {\bf 5} (2011), 22-36.

\bibitem{SA2} J.W. Sander and T. Sander,
{Integral circulant graphs of prime power order with maximal energy},
Linear Algebra Appl., to appear.    

\bibitem{SAX} N. Saxena and S. Severini and I.E. Shparlinski,
{Parameters of integral circulant graphs and periodic quantum dynamics},
Int. J. Quantum Inf. {\bf 5} (2007), 417-430.        

\bibitem{SHP} I. Shparlinski,
{On the energy of some circulant graphs},
Linear Algebra Appl. {\bf 414} (2006), 378-382.

\bibitem{SO} W. So,
Integral circulant graphs,
Discrete Math. {\bf 306} (2005), 153-158.

\bibitem{ST} K. Stolarsky, 
{Beatty sequences, continued fractions, and certain shift operators},
Canadian Math. Bull. {\bf 19} (1976), 473-482.

\end{thebibliography}
\end{document}